\documentclass{elsart}

\usepackage{amssymb}
\usepackage[latin1]{inputenc}

\begin{document}

\begin{frontmatter}



\title{Accurate numerical linear algebra with Bernstein--Vandermonde matrices}


\author{Ana Marco
\corauthref{1}}
\corauth[1]{Corresponding author.}
\ead{ana.marco@uah.es},
\author{José-Javier Mart{\'\i}nez}
\ead{jjavier.martinez@uah.es}

\address{Departamento de Matemáticas, Universidad de Alcalá, Campus Universitario, 28871 Alcalá de Henares, Madrid, Spain}

\begin{abstract}
The accurate solution of some of the main problems in numerical linear algebra (linear system solving, eigenvalue computation,  singular value computation and the least squares problem) for a totally positive Bernstein-Vandermonde matrix is considered. Bernstein-Vandermonde matrices are a generalization of Vandermonde matrices arising when considering for the space of the algebraic polynomials of degree less than or equal to $n$ the Bernstein basis, a widely used basis in Computer Aided Geometric Design, instead of the monomial basis.

Our approach is based on the computation of the bidiagonal factorization of a totally positive Bernstein-Vandermonde matrix (or its inverse) by means of Neville elimination. The explicit expressions obtained for the determinants involved in the process makes the algorithm  both fast and accurate.


\end{abstract}

\end{frontmatter}


\section{Introduction}

The design of accurate and efficient algorithms for structured matrices is a relevant field in numerical linear algebra which in recent years has received a growing attention (see, for instance, the recent survey paper \cite{DDHK} and references therein). An important condition for an algorithm to be accurate is the so called NIC (no inaccurate cancellation) condition, which is emphasized in \cite{DDHK}:

{\bf NIC:} The algorithm only multiplies, divides, adds (resp., substracts) real numbers with like (resp., differing) signs, and otherwise only adds or substracts input data.

In particular, in Section 2.3 of \cite{DDHK} several different classes of structured matrices are considered, among them the class of totally positive matrices. In our work we will consider one special class of totally positive matrices which we have recently analyzed in \cite{MM07} in the context of linear system solving: Bernstein-Vandermonde matrices.

Bernstein-Vandermonde matrices are a generalization of Vandermonde matrices arising when considering for the space of the algebraic polynomials of degree less than or equal to $n$ the {\it Bernstein basis}, a widely used basis in Computer Aided Geometric Design due to the good properties that it possesses (see, for instance, \cite{CP, FARIN}). The explicit
conversion between the Bernstein basis and the power basis is
exponentially ill-conditioned as the polynomial degree increases
\cite{FAROU}, and so it is very important that when designing
algorithms for performing numerical computations with polynomials
expressed in Bernstein form, all the intermediate operations are
developed using this form only.

As recalled in \cite{DDHK}, a crucial preliminary stage of the algorithms considered there for the class of totally positive matrices is the decomposition of the matrix as a product of bidiagonal factors, and that is our main contribution for the case of Bernstein-Vandermonde matrices: the fast and accurate computation of this {\it bidiagonal factorization}.

Working along this line, in this work we present algorithms for linear system solving, eigenvalue computation, singular value computation and the least squares problem which are based on the bidiagonal factorization of the
 Bernstein-Vandermonde matrix or of its inverse. Factorizations in terms
of bidiagonal matrices are very useful when working with Vandermonde
\cite{BP, HIG}, Cauchy \cite{BKO}, Cauchy-Vandermonde \cite{MP98, MP03} and generalized
Vandermonde matrices \cite{DK}.

As for the other stages of our algorithms, we will use several algorithms developed by Plamen Koev \cite{KOEV, KOEV05, KOEV07}.

The rest of the paper is organized as follows. Some basic results on Neville elimination
and total positivity are recalled in Section 2. In Section 3 the
bidiagonal factorization of a Bernstein-Vandermonde
matrix and of its inverse are considered. In Section 4 the algorithms for solving the four above problems are presented. Finally, Section 5 is
devoted to illustrate the accuracy of the algorithms by means of some
numerical experiments.

\section{Basic facts on Neville elimination and total positivity}

To make this paper as self-contained as possible, we will briefly recall in this section some basic results on Neville
elimination and total positivity which will be essential for obtaining
the results presented in Section 3. Our notation follows the notation used in \cite{GP92} and \cite{GP94}. Given $k$,
$n \in {\bf N}$  ($1 \leq k \leq n$), $Q_{k,n}$ will denote the set
of all increasing sequences of $k$ positive integers less than or
equal to $n$.

Let $A$ be an $l \times n$ real matrix. For $k \leq l$, $m
\leq n$, and for any $\alpha \in Q_{k,l}$ and $\beta \in Q_{m,n}$,
we will denote by $A[\alpha \vert \beta ]$ the submatrix $k\times m$
of $A$ containing the rows numbered by $\alpha $ and the columns
numbered by $\beta $.

The fundamental tool for obtaining the results presented in this
paper is the {\it Neville elimination} \cite{GP92,GP94}, a procedure that
makes zeros in a matrix adding to a given row an appropriate
multiple of the previous one. We will describe the Neville elimination for a matrix $A=(a_{i,j})_{1\leq i\leq l; 1\leq j\leq n}$ where $l \geq n$.

Let $A=(a_{i,j})_{1\leq i\leq l; 1\leq j\leq n}$ be a matrix where $l \geq n$. The Neville elimination of $A$ consists of $n-1$ steps
resulting in a sequence of matrices $A:=A_1\to A_2\to \ldots \to
A_n$, where $A_t= (a_{i,j}^{(t)})_{1\leq i\leq l; 1\leq j\leq n}$ has zeros below
its main diagonal in the $t-1$ first columns. The matrix $A_{t+1}$
is obtained from $A_t$ ($t=1,\ldots ,n$) by using the following
formula:
$$
a_{i,j}^{(t+1)}:= \left\{ \begin{array}{ll} a_{i,j}^{(t)}~, &
\textnormal{if} \quad i\leq t\\
a_{i,j}^{(t)}-(a_{i,t}^{(t)}/a_{i-1,t}^{t})a_{i-1,j}^{(t)}~,~ &
\textnormal{if} \quad i\geq t+1 ~\textnormal{and}~ j\geq t+1\\
0~, & \textnormal{otherwise}.
\end{array}\right.  \eqno(2.1)
$$
In this process the element
$$
p_{i,j}:=a_{i,j}^{(j)} \qquad 1\leq j\leq n, ~~ j\leq i\leq l
$$
is called {\it pivot} ($i,j$) of the Neville elimination of $A$. The
process would break down if any of the pivots $p_{i,j}$ ($1\leq j\leq n, ~~ j\leq i\leq l$) is zero. In that case we can move the corresponding rows to
the bottom and proceed with the new matrix, as described in \cite{GP92}.
The Neville elimination can be done without row exchanges if all the
pivots are nonzero, as it will happen in our situation. The pivots
$p_{i,i}$ are called {\it diagonal pivots}. If all the pivots
$p_{i,j}$ are nonzero, then $p_{i,1}=a_{i,1}\, \forall i$ and, by
Lemma 2.6 of \cite{GP92}
$$p_{i,j}={\det A[i-j+1,\ldots ,i\vert 1,\ldots ,j]\over \det
A[i-j+1,\ldots ,i-1\vert 1,\ldots ,j-1]}\qquad 1<j\leq n, ~ j\leq i\leq l.
\eqno(2.2)$$
The element
$$
m_{i,j}=\frac{p_{i,j}}{p_{i-1,j}} \qquad 1\leq j\leq n, ~~ j< i\leq
l  \eqno(2.3)
$$
is called {\it multiplier} of the Neville elimination of $A$. The
matrix $U:=A_n$ is upper triangular and has the diagonal pivots in
its main diagonal.

The {\it complete Neville elimination} of a matrix $A$ consists on
performing the Neville elimination of $A$ for obtaining $U$ and then
continue with the Neville elimination of $U^T$. The pivot
(respectively, multiplier) $(i,j)$ of the complete Neville
elimination of $A$ is the pivot (respectively, multiplier) $(j,i)$
of the Neville elimination of $U^T$, if $j\ge i$. When no row
exchanges are needed in the Neville elimination of $A$ and $U^T$, we
say that the complete Neville elimination of $A$ can be done without
row and column exchanges, and in this case the multipliers of the
complete Neville elimination of $A$ are the multipliers of the
Neville elimination of $A$ if $i\ge j$ and the multipliers of the
Neville elimination of $A^T$ if $j\ge i$.

A matrix is called {\it totally positive} (respectively, {\it
strictly totally positive}) if all its minors are nonnegative
(respectively, positive). The Neville elimination characterizes the
strictly totally positive matrices as follows \cite{GP92}:

{\bf Theorem 2.1.} A matrix is strictly totally positive if and only
if its complete Neville elimination can be performed without row and
column exchanges, the multipliers of the Neville elimination of $A$
and $A^T$ are positive, and the diagonal pivots of the Neville
elimination of $A$ are positive.

It is well known \cite{CP} that the Bernstein-Vandermonde matrix is a
strictly totally positive matrix when the nodes
satisfy $0<x_1<x_2<\ldots<x_{l+1}<1$,  but this result is also a consequence of our Theorem 3.2. \cite{MM08}.

\section{Bidiagonal factorizations} 

The {\it Bernstein basis} of the space $\Pi_n(x)$ of polynomials of
degree less than or equal to $n$ on the interval $[0,1]$ is:
$$
\mathcal{B}_n=\big\{ b_i^{(n)}(x) = {n \choose i} (1 - x)^{n-i} x^i,
\qquad i = 0, \ldots, n \big\}.
$$
The matrix
$$
A=
\left(\begin{array}{cccc}
{n \choose 0}(1-x_1)^n & {n \choose 1}x_1(1-x_1)^{n-1}& \cdots & {n
\choose n}x_1^n\\
{n \choose 0}(1-x_2)^n & {n \choose
1}x_2(1-x_2)^{n-1}& \cdots & {n \choose n}x_2^n\\
\vdots & \vdots & \ddots & \vdots\\ {n \choose 0}(1-x_{l+1})^n & {n
\choose 1}x_{l+1}(1-x_{l+1})^{n-1}& \cdots & {n \choose n}x_{l+1}^n
\end{array}\right)
$$
is the $(l+1)\times(n+1)$ {\it Bernstein-Vandermonde matrix} for the Bernstein basis $\mathcal{B}_n$ and the nodes $\{x_i \}_{1\leq i \leq l+1}$.

From now on, we will assume $0<x_1<x_2< \ldots <x_{l+1}<1$. In this situation the Bernstein-Vandermonde matrix is a strictly totally positive matrix \cite{CP} and the following two theorems hold:

{\bf Theorem 3.1.} Let $A=(a_{i,j})_{1\le i,j\le n+1}$ be a
Bernstein-Vandermonde matrix whose nodes satisfy $0 < x_1 < x_2
<\ldots < x_n <x_{n+1} <1$. Then $A^{-1}$ admits a factorization in
the form
$$A^{-1}=G_1G_2\cdots G_{n}D^{-1}F_{n}F_{n-1}\cdots F_1, \eqno(3.1)$$
where
$F_i$ ($i=1,\ldots,n$) are $(n+1)\times(n+1)$ bidiagonal matrices of the form
$$F_i=\left( \begin{array}{cccccccc}
1 & & & & & & & \\
0 & 1 & & & & & & \\
& \ddots & \ddots & & & & & \\
& & 0 & 1 & & & & \\
& & & -m_{i+1,i} & 1 & & & \\
& & & & -m_{i+2,i} & 1 & & \\
& & & & & \ddots & \ddots & \\
& & & & & & -m_{n+1,i} & 1
\end{array} \right),  \eqno(3.2)
$$
$G^T_i$ ($1\le i\le n$) are $(n+1)\times(n+1)$ bidiagonal matrices of the form
$$G_i^T=\left( \begin{array}{cccccccc}
1 & & & & & & & \\
0 & 1 & & & & & & \\
& \ddots & \ddots & & & & & \\
& & 0 & 1 & & & & \\
& & & -\widetilde m_{i+1,i} & 1 & & & \\
& & & & -\widetilde m_{i+2,i} & 1 & & \\
& & & & & \ddots & \ddots & \\
& & & & & & -\widetilde m_{n+1,i} & 1
\end{array} \right),  \eqno(3.3)
$$
$(i=1,\ldots,n)$, and $D$ is a diagonal matrix of order $n+1$
$$
D=\textnormal{diag}\{p_{1,1},p_{2,2},\ldots,p_{n+1,n+1}\}. \eqno(3.4)
$$
$m_{i,j}$ are the multipliers of the Neville elimination of the Bernstein-Vandermonde matrix $A$, and have the expression
$$
m_{i,j}=\frac{ (1-x_i)^{n-j+1} (1-x_{i-j})
\prod_{k=1}^{j-1}(x_i-x_{i-k}) }{ (1-x_{i-1})^{n-j+2}
\prod_{k=2}^j(x_{i-1}-x_{i-k}) },  \eqno(3.5)
$$
where $j=1,\ldots,n$ and $i=j+1, \dots, n+1$.

$\widetilde m_{i,j}$ are the multipliers of the Neville elimination of $A^T$ and their expression is
$$
\widetilde m_{i,j}=\frac{(n-i+2)\cdot x_j}{(i-1)(1-x_j)}, \eqno(3.6)
$$
where $j=1,\dots,n$ and $i=j+1,\ldots,n+1$.

Finally, the $i$th diagonal element of $D$ is the diagonal pivot of the Neville elimination of $A$ and its expression is
$$
p_{i,i}=\frac{ {n \choose i-1}(1-x_i)^{n-i+1} \prod_{k<i}(x_i-x_k)
}{ \prod_{k=1}^{i-1}(1-x_k) }  \eqno(3.7)
$$
for $i=1,\ldots,n+1$.

{\bf Proof.} It can be found in \cite{MM07}. $\Box$

{\bf Theorem 3.2.} Let $A=(a_{i,j})_{1\le i \leq l+1; 1\leq j\le n+1}$ be a
Bernstein-Vandermonde matrix for the Bernstein basis $\mathcal{B}_n$ whose nodes satisfy $0 < x_1 < x_2
<\ldots < x_l <x_{l+1} <1$. Then $A$ admits a factorization in
the form
$$A=F_{l}F_{l-1} \cdots F_1 D G_1 \cdots G_{n-1}G_{n}  \eqno(3.8)$$
where $F_i$ ($1\le i\le l$) are $(l+1)\times(l+1)$  bidiagonal matrices of the form
$$F_i=\left( \begin{array}{cccccccc}
1 & & & & & & & \\
0 & 1 & & & & & & \\
& \ddots & \ddots & & & & & \\
& & 0 & 1 & & & & \\
& & & m_{i+1,1} & 1 & & & \\
& & & & m_{i+2,2} & 1 & & \\
& & & & & \ddots & \ddots & \\
& & & & & & m_{l,l-i} & 1
\end{array}\right),  \eqno(3.9)
$$
$G^T_i$ ($1\le i\le n$) are $(n+1)\times (n+1)$ bidiagonal matrices of the form
$$G_i^T=\left(\begin{array}{cccccccc}
1 & & & & & & & \\
0 & 1 & & & & & & \\
& \ddots & \ddots & & & & & \\
& & 0 & 1 & & & & \\
& & & \widetilde m_{i+1,1} & 1 & & & \\
& & & & \widetilde m_{i+2,2} & 1 & & \\
& & & & & \ddots & \ddots & \\
& & & & & & \widetilde m_{n,n-i} & 1
\end{array}\right),  \eqno(3.10)
$$
and $D$ is the $(l+1)\times (n+1)$ diagonal matrix
$$
D =(d_{i,j})_{1\le i \leq l+1; 1\leq j\le n+1}=\textnormal{diag}\{p_{1,1},p_{2,2},\ldots,p_{n+1,n+1}\}. \eqno(3.11)
$$
The expressions of the multipliers $m_{i,j}$ $(j=1,\ldots,n+1;\quad i=j+1,\ldots,l+1)$ of the Neville elimination of $A$, the multipliers $\widetilde m_{i,j}$ $(j=1,\ldots,n;\quad i=j+1,\ldots,n+1)$ of the Neville elimination of $A^T$, and the diagonal pivots $P_{i,i}$ $(i=1,\ldots,n+1)$ of the Neville elimination of $A$ are also in this case the given by Eq. (3.5), Eq. (3.6) and Eq. (3.7), respectively.

{\bf Proof.} It can be found in \cite{MM08}. $\Box$

It must be observed that in the square case, the matrices $F_i$ ($i=1,\ldots, l$) and the matrices $G_j$ ($j=1,\ldots,n$) that appear in the bidiagonal factorization of $A$ are not the same bidiagonal matrices that appear in the bidiagonal factorization of $A^{-1}$ , nor their inverses (see Theorem 3.1 and Theorem 3.2). The multipliers of the Neville elimination of $A$ and $A^T$ give us the bidiagonal factorization of $A$ and $A^{-1}$, but obtaining the bidiagonal factorization of $A$ from the bidiagonal factorization of $A^{-1}$ (or vice versa) is not straightforward. The structure of the bidiagonal matrices that appear in both factorizations is not preserved by the inversion, that is, in general, $F_i^{-1}$ ($i=1,\ldots,l$) and $G_j^{-1}$ ($j=1,\ldots,n$) are not bidiagonal matrices. See \cite{GP96} for a more detailed explanation.

A fast and accurate algorithm for computing the bidiagonal factorization of the totally positive Bernstein-Vandermonde matrix $A$ and of its inverse (when it exists) has been developed by using the expressions (3.5), (3.6) and (3.7) for the computation of the multipliers $m_{i,j}$ and  $\widetilde m_{i,j}$, and the diagonal pivots $p_{i,i}$ of its Neville elimination \cite{MM07, MM08}. Given the nodes $\{x_i \}_{1\leq i \leq l+1} \in (0,1)$ and the degree $n$ of the Bernstein basis, it returns a matrix $M \in {\bf R}^{(l+1)\times(n+1)}$ such that
$$
\begin{array}{l}
M_{i,i} =p_{i,i} \quad ~~i=1,\ldots,n+1,\\
M_{i,j} =m_{i,j} \quad  j=1,\ldots, n+1; ~i=j+1,\ldots, l+1,\\
M_{i,j} =\widetilde m_{j,i} \quad i=1,\ldots,n; ~j=i+1,\ldots, n+1.
\end{array}
$$
The algorithm, which we have called {\tt TNBDBV}, does not construct the Bernstein-Vandermonde matrix, it only works with the nodes $\{x_i\}_{1\leq i\leq l+1}$. Its computational cost is of $O(ln)$ arithmetic operations, and has high relative accuracy because it only involves arithmetic operations that avoid inaccurate cancellation (see \cite{MM07} for the details). The implementation in \textsc{Matlab} of the algorithm in the square case can be taken from \cite{KOEV}.

{\it Remark.} The algorithm {\tt TNBDBV} computes the matrix $M$, denoted
as $\mathcal{BD}(A)$ in \cite{KOEV05}, which represents the {\it bidiagonal
decomposition} of $A$. But it is a remarkable fact that the same
matrix $\mathcal{BD}(A)$ also serves to represent the bidiagonal
decomposition of $A^{-1}$.

\section{Accurate computations with Bernstein-Vandermonde matrices} 

In this section four fundamental problems in numerical linear algebra (linear system solving, eigenvalue computation, singular value computation and the least squares problem) are considered for the case of a totally positive Bernstein-Vandermonde matrix. The bidiagonal factorization of the Bernstein-Vandermonde matrix (or its inverse) let us to develop accurate and efficient algorithms for solving each one of these problems.

Let us observe here that, of course, one could try to solve these problems by using standard algorithms. However 
 the solution provided by them will generally be less accurate since
 Bernstein-Vandermonde matrices are ill conditioned (see \cite{MM07}) and
these algorithms can suffer from inaccurate cancellation, since they do not take into account the structure of
the matrix, which is crucial in our approach.

\subsection{Linear system solving}

Let $Ax=b$ be a linear system whose coefficient matrix $A$ is a square Bernstein-Vandermonde matrix of order $n+1$ generated by the nodes $\{x_i\}_{1\leq i\leq n+1}$, where $0<x_1<\ldots<x_{n+1}<1$. An application which involves the solution of this type of linear systems has been presented in \cite{MM07a}.

The following algorithm solves $Ax=b$ accurately with a computational cost of $O(n^2)$ arithmetic operations (see \cite{MM07} for the details):

INPUT: The nodes $\{x_i\}_{1\leq i\leq n+1}$ and the data vector $b\in {\bf R}^{n+1}$.

OUTPUT: The solution vector $x\in {\bf R}^{n+1}$.

\begin{itemize}
\item[-] {\it Step 1:} Computation of the bidiagonal decomposition of $A^{-1}$ by using {\tt TNBDBV}.

\item[-] {\it Step 2:} Computation of $$x=A^{-1}b=G_1G_2\cdots G_{n}D^{-1}F_{n}F_{n-1}\cdots F_1b$$
\end{itemize}

Step 2 can be carried out by using the algorithm {\tt TNSolve} of P. Koev \cite{KOEV}. Given the bidiagonal factorization of the matrix $A$, {\tt TNSolve} solves $Ax=b$ accurately by using backward substitution.

Several examples illustrating the good behaviour of our algorithm can be found in \cite{MM07}.

\subsection{Eigenvalue computation}

Let $A$ be a square Bernstein-Vandermonde matrix of order $n+1$ generated by the nodes $\{x_i\}_{1\leq i\leq n+1}$, where $0<x_1<\ldots<x_{n+1}<1$. The following algorithm computes accurately the eigenvalues of $A$.

INPUT: The nodes $\{x_i\}_{1\leq i\leq n+1}$.

OUTPUT: A vector $x\in {\bf R}^{n+1}$ containing the eigenvalues of $A$.

\begin{itemize}
\item[-] {\it Step 1:} Computation of the bidiagonal decomposition of $A^{-1}$ by using {\tt TNBDBV}.

\item[-] {\it Step 2:} Given the result of Step 1, computation of the eigenvalues of $A$ by using the algorithm {\tt TNEigenvalues}.
\end{itemize}

{\tt TNEigenvalues} is an algorithm of P. Koev \cite{KOEV05} which computes accurate eigenvalues of a totally positive matrix starting from its bidiagonal factorization. The computational cost of {\tt TNEigenvalues} is of $O(n^3)$ arithmetic operations (see \cite{KOEV05}) and its implementation in \textsc{Matlab} can be taken from \cite{KOEV}. In this way, as the computational cost of Step 1 is of $O(n^2)$ arithmetic operations, the cost of the whole algorithm is of $O(n^3)$ arithmetic operations.

\subsection{The least squares problem}

Let $A \in {\bf R}^{(l+1)\times(n+1)}$ be a Bernstein-Vandermonde matrix generated by the nodes $\{x_i\}_{1\leq i\leq l+1}$, where $0<x_1<\ldots<x_{l+1}<1$ and $l>n$. Let $b\in{\bf R}^{l+1}$ be a data vector. The least squares problem associated to $A$ and $b$ consists on  computing a vector $x \in {\bf R}^{n+1}$ minimizing $\parallel Ax-b \parallel_2$.

Taking into account that in the situation we are considering $A$ is a strictly totally positive matrix, it has full rank, and the method based on the QR decomposition due to Golub \cite{GOL} is adequate \cite{BJO}. For the sake of completeness, we include the following result (see Section 1.3.1 in \cite{BJO}) which will be essential in the construction of our algorithm.

{\bf Theorem 4.1.} Let $Ac=f$ a linear system where $A \in {\bf R}^{(l+1) \times (n+1)}$, $l \geq n$, $c \in {\bf R}^{n+1}$ and $f \in {\bf R}^{l+1}$. Assume that $rank(A)=n+1$, and let the QR decomposition of $A$  be given by
$$
A=Q \left(\begin{array}{c} R\\ 0 \end{array}\right),
$$
where $Q \in {\bf R}^{(l+1) \times (l+1)}$ is an orthogonal matrix and $R \in {\bf R}^{(n+1) \times (n+1)}$ is an upper triangular matrix with positive diagonal entries. Then the solution of the least squares problem $min_c \parallel f-Ac \parallel_2$ is obtained from
$$
\left(\begin{array}{c} d_1 \\ d_2\end{array}\right) = Q^T f, \quad Rc=d_1, \quad r=Q \left(\begin{array}{c}0 \\ d_2 \end{array}\right),
$$
where $d_1 \in {\bf R}^{n+1}$, $d_2\in {\bf R}^{l-n}$ and $r=f-Ac$. In particular $\parallel r \parallel_2 = \parallel d_2 \parallel_2$.

The following algorithm, which is based on the previous theorem, solves in accurate and efficient way our least squares problem:

INPUT: The nodes $\{x_i \}_{1\leq i \leq l+1}$, the data vector $f$ and the degree $n$ of the Bernstein basis.

OUTPUT: The vector $x \in {\bf R}^{n+1}$ minimizing $\parallel Ax-b \parallel_2$ and the minimun residual $r=b-Ax$.

\begin{itemize}

\item[-] {\it Step 1}: Computation of the bidiagonal factorization of $A$ by means of {\tt TNBDBV}.

\item[-] {\it Step 2}: Given the result of Step 1, computation of the QR decomposition of $A$ by using {\tt TNQR}.

\item[-] {\it Step 3}: Computation of
$$d=\left( \begin{array}{c}d_1 \\d_2 \end{array} \right)=Q^Tf.$$

\item[-] {\it Step 4}: Solution of the upper triangular system $Rc=d_1$.

\item[-] {\it Step 5}: Computation of
$$
r=Q \left(\begin{array}{c}0 \\ d_2 \end{array}\right).
$$

\end{itemize}

The algorithm {\tt TNQR} has been developed by P. Koev, and given the bidiagonal factorization of $A$, it computes the matrix $Q$ and the bidiagonal factorization of the matrix $R$. Let us point out here that if $A$ is strictly totally positive, then $R$ is strictly totally positive. {\tt TNQR} is based on Givens rotations, has a computational cost of $O(l^2n)$ arithmetic operations if the matrix $Q$ is required, and its high relative accuracy comes from the avoidance of inaccurate cancellation \cite{KOEV07}. Its implementation in \textsc{Matlab} can be obtained from \cite{KOEV}.

As for the computational cost of the whole algorithm, it is led by the cost of computing the QR decomposition of $A$, and therefore, it is of $O(l^2n)$ arithmetic operations.

Some numerical experiments which show the good behaviour of our algorithm when solving problems of polynomial regression in the Bernstein basis have been presented in \cite{MM08}.

\subsection{Singular value computation}

Let $A\in {\bf R}^{(l+1)\times(n+1)}$ be a Bernstein-Vandermonde matrix generated by the nodes $\{x_i\}_{1\leq i\leq l+1}$, where $0<x_1<\ldots<x_{l+1}<1$ and $l>n$. The following algorithm computes in an accurate and efficient way the singular values of $A$.

INPUT: The nodes $\{x_i\}_{1\leq i\leq l+1}$ and the degree $n$ of the Bernstein basis.

OUTPUT: A vector $x\in{\bf R}^{n+1}$ containing the singular values of $A$.

\begin{itemize}
\item[-]{\it Step 1:} Computation of the bidiagonal decomposition of $A$ by using {\tt TNBDBV}.

\item[-]{\it Step 2:} Given the result of Step 1, computation of the singular values by using {\tt TNSingularValues}.
\end{itemize}

{\tt TNSingularValues} is an algorithm of P. Koev that computes accurate singular values of a totally positive matrix starting from its bidiagonal factorization \cite{KOEV05}.  Its computational cost is of $O(ln^2)$ and its implementation in \textsc{Matlab} can be found in \cite{KOEV}. Taking this complexity into account, the computational cost of our algorithm for computing the singular values of a totally positive Bernstein-Vandermonde matrix is of $O(ln^2)$ arithmetic operations.

\section{Numerical experiments}

In this last section we include two numerical experiments which illustrate the high relative accuracy of the algorithms we have presented for the problems of eigenvalue computation and singular value computation. Numerical experiments for the cases of linear system solving and of least squares problems have been included in \cite{MM07, MM08}. 

{\bf Example 5.1.} Let ${\mathcal B}_{20}$ be the Bernstein basis of the space of polynomials with degree less than or equal to $20$ on $[0, 1]$ and let$A$ be the square Bernstein-Vandermonde matrix of order $21$ generated by the nodes:

$\frac{1}{22}< \frac{1}{20}< \frac{1}{18}< \frac{1}{16}< \frac{1}{14}< \frac{1}{12}< \frac{1}{10}< \frac{1}{8}< \frac{1}{6}< \frac{1}{4}< \frac{1}{2}<
\frac{23}{42}< \frac{21}{38}< \frac{19}{34}< \frac{17}{30}< \frac{15}{26}< \frac{13}{22}< \frac{11}{18}< \frac{9}{14}< \frac{7}{10}< \frac{5}{6}.$

The condition number of $A$ is $\kappa_2(A)=1.9e+12$. In Table 1 we present the eigenvalues $\lambda_i$ of $A$ and the relative errors obtained when computing them by means of:
\begin{enumerate}
\item The algorithm presented in Section 4.2 (column labeled by {\tt MM}).
\item The command {\tt eig} from \textsc{Matlab}.
\end{enumerate}
The relative error of each computed eigenvalue is obtained by using the eigenvalues calculated in {\it Maple 10} with 50-digit arithmetic.

\begin{table}[h]
\begin{center}
\begin{tabular}{|c|c|c|}
\hline $\lambda_i$  & {\tt MM} & {\tt eig} \\
\hline $1.0e+00$ & $0$ & $4.0e-15$ \\
\hline $8.4e-01$ & $4.0e-16$ & $1.3e-16$ \\
\hline $2.8e-01$ & $9.9e-16$ & $3.0e-15$ \\
\hline $2.1e-01$ & $0$ & $3.7e-15$ \\
\hline $1.2e-01$ & $6.0e-16$ & $3.6e-15$ \\
\hline $6.6e-02$ & $8.5e-16$ & $5.5e-15$ \\
\hline $3.8e-02$ & $3.7e-16$ & $1.2e-14$ \\
\hline $2.2e-02$ & $4.8e-16$ & $2.0e-14$ \\
\hline $9.4e-03$ & $0$ & $4.6e-14$ \\
\hline $4.6e-03$ & $1.9e-16$ & $9.3e-14$ \\
\hline $1.5e-03$ & $2.9e-16$ & $6.1e-14$ \\
\hline $5.9e-04$ & $1.8e-16$ & $2.5e-13$ \\
\hline $1.7e-04$ & $3.2e-16$ &  $1.3e-12$ \\
\hline $4.3e-05$ & $9.4e-16$ &  $3.6e-12$ \\
\hline $1.3e-05$ & $1.1e-15$ & $6.9e-12$ \\
\hline $1.7e-06$ & $6.1e-16$ & $5.1e-11$\\
\hline $5.6e-07$ & $0$ & $1.9e-10$ \\
\hline $3.5e-08$ & $2.8e-15$ & $1.8e-09$ \\
\hline $1.1e-08$ & $5.9e-16$ & $5.9e-09$ \\
\hline $2.7e-10$ & $2.1e-15$ & $6.4e-08$\\
\hline $1.3e-12$ & $9.0e-16$ & $1.0e-05$\\
\hline
\end{tabular}
\end{center}\caption{Example 5.1: Eigenvalues of a Bernstein-Vandermonde matrix of order $21$.}
\end{table}

{\bf Example 5.2.} Let ${\mathcal B}_{15}$ be the Bernstein basis of the space of polynomials with degree less than or equal to $16$ on $[0, 1]$ and let $A\in {\bf}^{21\times 16}$ Bernstein-Vandermonde matrix generated by the same nodes as in Example 5.1. Its condition number is $\kappa_2(A)=5.3e+08$. In Table 2 we present the singular values $\sigma_i$ of $A$ and the relative errors obtained when computing them by means of
\begin{enumerate}
\item The algorithm presented in Section 4.4 (column labeled by {\tt MM}).
\item The command {\tt svd} from \textsc{Matlab}.
\end{enumerate}
The relative error of each computed singular value is obtained by using the singular values calculated in {\it Maple 10} with 50-digit arithmetic.

\begin{table}[h]
\begin{center}
\begin{tabular}{|c|c|c|}
\hline $\sigma_i$  & {\tt MM} & {\tt svd} \\
\hline $1.6e+00$ & $1.4e-16$ & $1.4e-16$ \\
\hline $1.2e+00$ & $9.4e-16$ & $9.4e-16$ \\
\hline $5.8e-01$ & $5.8e-16$ &  $3.8e-16$ \\
\hline $4.8e-01$ & $5.8e-16$ &  $2.3e-16$ \\
\hline $2.5e-01$ & $2.2e-16$ &  $0$ \\
\hline $1.8e-01$ & $6.0e-16$ &  $4.5e-16$ \\
\hline $6.3e-02$ & $2.2e-16$ &  $0$ \\
\hline $3.5e-02$ & $0$ &  $1.0e-15$ \\
\hline $7.6e-03$ & $1.1e-16$ & $9.8e-15$ \\
\hline $3.4e-03$ & $1.1e-16$ & $2.0e-15$ \\
\hline $3.9e-04$ & $1.4e-16$ & $1.5e-14$ \\
\hline $1.3e-04$ & $2.0e-16$ & $8.4e-15$ \\
\hline $7.5e-06$ & $9.0e-16$ & $1.6e-12$ \\
\hline $1.6e-06$ & $1.32e-16$ & $4.3e-12$ \\
\hline $4.8e-08$ & $1.8e-16$ &  $3.9e-10$ \\
\hline $3.0e-09$ & $2.9e-15$ &  $2.5e-10$ \\
\hline
\end{tabular}
\end{center}\caption{Example 5.2: Singular values of a Bernstein-Vandermonde matrix $21\times 16$.}
\end{table}

{\bf Acknowledgements.} This research has been partially supported by Spanish Research Grant
MTM2006-03388 from the Spanish Ministerio de Educaci\'on y Ciencia.

\label{}



\end{document}